\documentclass[12pt,draft]{amsart}

\usepackage{amssymb}

\begin{document}

\newtheoremstyle{hplain}%
  {\topsep}
  {\topsep}
  {\itshape}
  {}
  {\bfseries}
  {.}
  { }
  {\thmname{#2 }\thmnumber{#1}\thmnote{ \rm(#3)}}%
  
\newtheoremstyle{hdefinition}
  {\topsep}%
  {\topsep}%
  {\upshape}
  {}%
  {\bfseries}%
  {.}
  { }%
  {\thmname{#2 }\thmnumber{#1}\thmnote{ \rm(#3)}}%
 
\theoremstyle{hplain}
\newtheorem{thm}{Theorem}
\newtheorem{note}[thm]{Note}
\newtheorem{lemma}[thm]{Lemma}
\newtheorem{remark}[thm]{Remark}
\newtheorem{corollary}[thm]{Corollary}
\newtheorem{claim}{Claim}[thm]

\theoremstyle{hdefinition}
\newtheorem{definition}[thm]{Definition}
\newtheorem{question}[thm]{Question}

\newcommand{\restr}{\restriction}
\newcommand{\arr}{\longrightarrow}
\newcommand{\sub}{\subseteq}
\newcommand{\stick}{{}\sp {\bullet}\hskip -7.2pt \shortmid }
\newcommand{\squdia}{\square\hskip -13.5pt\diamondsuit}
\newcommand{\oom}{\stackrel{1-1}{\longrightarrow}}
\newcommand{\iso}{\stackrel{\sim}{\longrightarrow}}
\newcommand{\conc}{{}^\smallfrown}
\newcommand{\setm}{\setminus}
\newcommand{\Hx}[1]{{\rm H}_{\omega_{#1}}}
\newcommand{\Hk}[1]{{\rm H}_{#1}}
\newcommand{\Ht}{{\rm H}_\theta}
\newcommand{\Hl}{{\rm H}_\lambda}
\newcommand{\mbb}{\mathbb}
\newcommand{\mcal}{\mathcal}
\newcommand{\cc}{\twoheadrightarrow}
\newcommand{\notcc}{\twoheadrightarrow\!\!\!\!\!\!\!/\,\,\,\,\,\,\,}
\renewcommand{\frak}{\mathfrak}
\renewcommand{\theequation}{\thesection.\arabic{equation}}

\newcommand{\GSC}{\Gamma_\Sigma}

\numberwithin{equation}{section}

\newcommand{\rem}[1]{}
\newcommand{\comment}[1]{\marginpar{\scriptsize\raggedright #1}}
\renewcommand{\thesubsection}{\alph{subsection}}

\title[Forcing indestructibility]{Forcing indestructibility of
  set-theoretic axioms}
\author{Bernhard K\"onig}
\address{\newline
  Universit\'e Paris 7\newline
  2 place Jussieu\newline
  75251 Paris Cedex 05\newline
  France}

\thanks{The author acknowledges a grant awarded by the French
  ministry of research.}


\begin{abstract}
  Various theorems for the preservation of set-theoretic axioms under
  forcing are proved, regarding both forcing axioms and axioms true in
  the Levy-Collapse. These show in particular that certain
  applications of forcing axioms require to add generic countable
  sequences high up in the set-theoretic hierarchy even before
  collapsing everything down to $\aleph_1$. Later we give
  applications, among them the consistency of ${\rm MM}$ with
  $\aleph_\omega$ not being Jonsson which answers a question raised
  during Oberwolfach 2005.
\end{abstract}

\maketitle

\section{Introduction}

It was a widely held intuition in the early days of proper forcing
that a typical application of the Proper Forcing Axiom makes use of a
poset of the form $\sigma\mbox{-closed}\,*\,{\rm ccc}$. The usual
argument was to collapse the size of all objects to $\aleph_1$, then
use a ccc-poset to force the desired property to these objects and
finally pull everything back into the ground model with the help of
the forcing axiom. Later it was realized that forcing axioms have a
lot more applications than that. These new developments were heading
into two different directions, on the one hand there was the
development of semiproper forcing in \cite{shelah82:_proper_forcin}
which lead to the Semiproper Forcing Axiom and later to Martin's
Maximum in \cite{foreman88}. On the other hand, even for ${\rm PFA}$
it was soon realized that there is a large variety of proper forcings
that are not of the form $\sigma\mbox{-closed}\,*\,{\rm ccc}$.
Interesting examples here are posets that are not $\omega$-proper and
it was demonstrated first in \cite{shelah82:_proper_forcin} and later
in \cite{moore05:_set} that these can be used to good account. The
point we are trying to make is slightly different and includes both
proper and semiproper forcing notions. We give examples to show that
certain applications of forcing axioms can require adding reals or
countable sequences even before we collapse the size of the relevant
objects to $\aleph_1$. The intuition here is that if we collapse
without adding countable sequences then our object will have an
enumeration of order-type $\omega_1$ whose initial segments are in the
ground model.  But certain applications exclude such a possibility,
most prominently the negation of approachability properties. We just
sketched the general direction of this article. Section
\ref{sec:closure-prop} introduces a wide range of forcing properties
that will become important later for the indestructibility theorems.
In Section \ref{sec:indestruct} we will give various indestructibility
results for forcing axioms, but also indestructibility results for
axioms true in the Levy-Collapse of a large cardinal, most notably the
axiom of Game Reflection from \cite{koenig04:_gener}. Section
\ref{sec:Kurepa-tree} is a way of applying this technique to Jonsson
cardinals and related model-theoretic transfer properties. Among other
things, it is shown there that ${\rm MM}$ is consistent with
$\aleph_\omega$ not being Jonsson which shoots down a lingering
conjecture.

The reader is assumed to have a strong background in set theory. As a
general reference we recommend \cite{Kunen-intro} and as a reference
regarding proper and semiproper forcing we suggest
\cite{baumgartner84:_applic_proper_forcin_axiom} and
\cite{shelah98:_proper_improper}. The remaining paragraphs of the
introduction will be used to give a short summary of the most
important Lemmas and Definitions used in this paper.

\begin{definition}
  If $\Gamma$ is a class of posets then ${\rm MA}(\Gamma)$ denotes the
  statement that whenever $\mbb{P}\in\Gamma$ and
  $D_\xi\;(\xi<\omega_1)$ is a collection of dense subsets of
  $\mbb{P}$ then there exists a filter $G\sub\mbb{P}$ such that
  $D_\xi$ intersects $G$ for all $\xi<\omega_1$. In particular, {\rm
    PFA} is ${\rm MA}(\mbox{proper})$ and ${\rm MM}$ is ${\rm
    MA}$(preserving stationary subsets of $\omega_1$). The interested
  reader is referred to
  \cite{baumgartner84:_applic_proper_forcin_axiom} and
  \cite{foreman88} for more history on these {\em forcing axioms}.
\end{definition}

\begin{definition}\label{def-BMA}
  We would like to remind the reader of the notion of a {\em bounded
    forcing axiom}. Assume that $\lambda$ is a cardinal and $\Gamma$ a
  class of posets, then the following are equivalent (see
  \cite{bagaria00:_bound} and \cite{todorcevic02:_local_pfa}):
  \begin{enumerate}
  \item For every $\mbb{P}\in\Gamma$ and a collection
    $D_\xi\;(\xi<\omega_1)$ of dense subsets of $\mbb{P}$ of size
    $\leq\lambda$ there is a filter $G\sub\mbb{P}$ such that $D_\xi$
    intersects $G$ for all $\xi<\omega_1$.
  \item For every $A\sub\lambda$ and every $\Sigma_1$-formula
    $\varphi(x)$, if there is some $\mbb{P}\in\Gamma$ such that
    $\Vdash_\mbb{P}\varphi(A)$ then there are stationarily many
    $M\prec\Hk{\lambda^+}$ of size $\aleph_1$ containing $A$ such that
    $\Hk{\lambda^+}\models\varphi(\pi_M(A))$, where $\pi_M$ is the
    transitive collapse of $M$.
  \end{enumerate}
  So let us denote the equivalent statements (1) and (2) by ${\rm
    MA}(\Gamma,\lambda)$ and for simplicity we write ${\rm
    PFA}(\lambda)$ for ${\rm MA}(\mbox{proper},\lambda)$ and ${\rm
    MM}(\lambda)$ for ${\rm MA}$(preserving stationary subsets of
  $\omega_1,\lambda$). The axioms ${\rm PFA}(\omega_1)$ and ${\rm
    MM}(\omega_1)$ are often denoted by ${\rm BPFA}$ and ${\rm BMM}$
  respectively. If $\Gamma$ is any class of posets, we also write
  ${\rm PFA}(\Gamma)$ for ${\rm MA}(\mbox{proper and in $\Gamma$})$
  and similarly with ${\rm MM}(\Gamma)$.
\end{definition}

\begin{definition}
  The {\em approachability property for $\kappa$} (${\rm
    AP}_{\kappa}$) is the statement that there is a sequence
  $(C_\alpha:\alpha<\kappa^+)$ such that for any $\alpha<\kappa^+$:
  \begin{enumerate}
  \item[(a)] $C_\alpha\sub\kappa^+$, ${\rm otp}\;C_\alpha\leq\kappa$,
  \end{enumerate}
  and there is a club $C\sub\lim(\kappa^+)$ such that for every
  $\gamma\in C$:
  \begin{enumerate}
  \item[(b)] $C_\gamma\sub\gamma$ is club,
  \item[(c)] the initial segments of $C_\gamma$ are enumerated before
    $\gamma$,\\ i.e. $\forall\alpha<\gamma\;\,\exists\beta<\gamma\;\,
    C_\gamma\cap\alpha=C_\beta$.
  \end{enumerate}
  A straightforward argument shows that ${\rm AP}_\kappa$ follows from
  either $\square_\kappa$ or from the cardinal arithmetic
  $\kappa^{<\kappa}=\kappa$.
\end{definition}

We use an abbreviation in the context of elementary embeddings:
$j:M\arr N$ means that $j$ is a non-trivial elementary embedding from
$M$ into $N$ such that $M$ and $N$ are transitive. The {\em critical
  point} of such an embedding, i.e. the first ordinal moved by $j$, is
denoted by ${\rm cp}(j)$. We write $jx$ for $j(x)$ in a context where
too many parentheses might be confusing. Let us remind ourselves of
the well-known extension Lemma for elementary embeddings first noticed
by Silver:

\begin{lemma}[Extension Lemma]\label{silverlem}
  Let $j:M\arr N$ and assume that $G\sub\mathbb{P}$ is generic over
  $M$ and $K\sub j(\mathbb{P})$ generic over $N$. If $j''G\sub K$ then
  there is a unique extension $j^*:M[G]\arr N[K]$ of $j$ such that
  $j^*(G)=K$.
\end{lemma}
\begin{proof}
  For each $\mathbb{P}$-name $\dot{\tau}$ simply let
  $j^*(\dot{\tau}[G])=j(\dot{\tau})[K]$.
\end{proof}

\begin{definition}
  Recall that the model-theoretic transfer property
  $$(\lambda_1,\lambda_0)\cc(\kappa_1,\kappa_0)$$
  means that every
  structure $(\lambda_1,\lambda_0,f_i)_{i<\omega}$ has an elementary
  substructure of the form $(A_1,A_0,f_i)_{i<\omega}$, where
  $|A_1|=\kappa_1$ and $|A_2|=\kappa_0$. The relation
  $(\omega_2,\omega_1)\cc(\omega_1,\omega)$ is usually called {\em
    Chang's conjecture}. A cardinal $\mu$ is called {\em
    $\kappa$-Rowbottom} if for all $\lambda<\mu$ we have
  $(\mu,\lambda)\cc(\mu,<\kappa)$. A cardinal $\mu$ is called {\em
    Jonsson} if every algebra of size $\mu$ has a proper subalgebra of
  size $\mu$.
\end{definition}

We also need the following well-known Lemmas:

\begin{lemma}\label{coll-iso}
  Let $\lambda$ be regular and assume that $\mcal{P}$ is a
  $\sigma$-closed poset of size $\lambda^{\aleph_0}$ that collapses
  $\lambda^{\aleph_0}$ to $\aleph_1$. Then $\mcal{P}$ is
  forcing-isomorphic to ${\rm Col}(\aleph_1,\lambda)$.\hfill\qed
\end{lemma}

\begin{lemma}\label{str-embeds}
  Assume that $\mcal{P}$ is strategically $\sigma$-closed and
  $\lambda\geq|\mcal{P}|$. Then $\mcal{P}$ completely embeds into
  ${\rm Col}(\aleph_1,\lambda)$.
\end{lemma}
\begin{proof}
  Clearly, $\mcal{P}$ completely embeds into $\mcal{P}\times{\rm
    Col}(\aleph_1,\lambda)$. Results in
  \cite{foreman83:_games_boolean} imply that $\mcal{P}\times{\rm
    Col}(\aleph_1,\lambda)$ is $\sigma$-closed. Finally, Lemma
  \ref{coll-iso} concludes that $\mcal{P}\times{\rm
    Col}(\aleph_1,\lambda)\cong{\rm Col}(\aleph_1,\lambda)$ and we are
  done.
\end{proof}

\section{$\omega_2$-closure properties}
\label{sec:closure-prop}

In this section we introduce five different properties of forcings
which all entail that no new $\omega_1$-sequences be added. We give a
small overview before defining them one by one, the following list is
increasing in logical strength:

\begin{enumerate}
\item $\omega_2$-distributive
\item weakly $(\omega_1+1)$-game-closed
\item strongly $(\omega_1+1)$-game-closed
\item $\omega_2$-closed
\item $\omega_2$-directed-closed
\end{enumerate}

So let us start with distributivity. A poset $\mbb{P}$ is called {\em
  $\kappa$-distributive} if the intersection of less than $\kappa$
many dense open subsets is again dense open. Note that this is
equivalent to saying that $\mbb{P}$ adds no new sequences of length
less than $\kappa$. It will become clear later why we are mostly
interested in the case $\kappa=\omega_2$. We have the following
proposition:

\begin{lemma}\label{pres-bounded}
  Let $\lambda\geq\aleph_1$. The bounded forcing axiom ${\rm
    MM}(\lambda)$ is preserved by $\lambda^+$-distri\-butive forcings.
  Moreover, ${\rm PFA}(\lambda)$ is preserved by proper
  $\lambda^+$-distributive forcings.
\end{lemma}
\begin{proof}
  This follows simply from the fact that $\lambda^+$-distributive
  forcings add no new elements to $\Hk{\lambda^+}$, so check that (2)
  of Definition \ref{def-BMA} holds in any (proper)
  $\lambda^+$-distributive extension. For the ${\rm
    MM}(\lambda)$-argument, note that $\lambda^+$-distributive
  forcings preserve stationary subsets of $\omega_1$.
\end{proof}

Now recall longer versions of the Banach-Mazur game on a partial
ordering $\mbb{P}$: \medskip

\begin{tabular}{c|cccccc}
  Empty & $p_0$ & $p_2$ & $\ldots$ & $p_\xi$ & $\ldots$\\\hline
  Nonempty & $\qquad p_1$ & $\qquad p_3$ & $\qquad\ldots$ & $\qquad
  p_{\xi+1}$ & $\qquad\ldots$
\end{tabular}
  
\medskip\noindent where $p_\xi\;(\xi<\alpha)$ is descending in
$\mbb{P}$ and Nonempty wins the game of length $\alpha$ if he can play
$\alpha$ times.

\begin{definition}\label{def-game-closed}
  A poset $\mbb{P}$ is called {\em weakly $\alpha$-game-closed} if
  Player Non\-empty has a winning strategy in the Banach-Mazur game of
  length $\alpha$, where Nonempty is allowed to play at limit stages.
  $\mbb{P}$ is called {\em strongly $\alpha$-game-closed} if Player
  Nonempty has a winning strategy in the same game except where Empty
  is allowed to play at limit stages.
\end{definition}

It is clear that strongly $(\kappa+1)$-game-closed posets are also
weakly $(\kappa+1)$-game-closed. Remember that the standard forcing to
add a $\square_\kappa$-sequence with initial segments is weakly
$(\kappa+1)$-game-closed. An ${\rm AP}_\kappa$-sequence can be added
with a considerably milder forcing. The following crucial fact is
originally from \cite{yoshinubo03:_approac}:

\begin{lemma}\label{force-AP}
  For all cardinals $\kappa$ there is a strongly
  $(\kappa+1)$-game-closed forcing $\mbb{A}_\kappa$ that adds an ${\rm
    AP}_\kappa$-sequence.\hfill\qed
\end{lemma}

\begin{definition}\label{plays-def}
  Assume for the following that $\mbb{P}$ is strongly
  $(\omega_1+1)$-game-closed. Let us fix a winning strategy $\sigma$
  for Nonempty in the Banach-Mazur game on $\mbb{P}$. Instead of
  forcing with $\mbb{P}$, we could also add a play of the game
  generically. Then this play induces a generic filter for $\mbb{P}$.
  Define $$\mbb{R}=\{\langle
  p_\xi:\xi\leq\gamma\rangle:\gamma<\omega_1 \mbox{ and }
  p_\xi\;(\xi\leq\gamma)\mbox{ is a partial $\sigma$-play}\}.$$
  If
  $s=\langle p_\xi:\xi\leq\gamma\rangle\in\mbb{R}$ is such a partial
  play, we also denote the maximal condition $p_\gamma$ by $p_s$. The
  ordering on $\mbb{R}$ is usual extension. Note that $\mbb{R}$ is
  $\sigma$-closed and contains $\mbb{P}$ as a complete subalgebra by
  the projection mapping $i(s)=p_s$. Yet, it is a much stronger
  forcing: $\mbb{R}$ will typically collapse the cardinality of
  $\mbb{P}$ to $\aleph_1$.
\end{definition}

\begin{lemma}\label{finsegsigma}
  Using the notation of Definition \ref{plays-def}, if $G\sub\mbb{P}$
  is generic then $\mbb{R}/\mbb{P}=\{s\in\mbb{R}:p_s\in G\}$ is
  $\sigma$-closed.
\end{lemma}
\begin{proof}
  Suppose $s_n\;(n<\omega)$ is a descending sequence in
  $\mbb{R}/\mbb{P}$ and $\gamma$ the length of the union
  $\bigcup_{n<\omega}s_n$. Then $q=\inf_{n<\omega}p_{s_n}$ is in $G$
  and $$s=\bigcup_{n<\omega} s_n\cup\{(\gamma,q)\}$$
  is still a
  partial play according to $\sigma$.
\end{proof}

\begin{lemma}\label{quot-sigma}
  Assume $\mbb{P}$ is strongly $(\omega_1+1)$-game-closed and
  $\lambda\geq|\mbb{P}|^{\aleph_0}$. Then ${\rm
    Col}(\omega_1,\lambda)/\mbb{P}$ is $\sigma$-closed.
\end{lemma}
\begin{proof}
  Let $\mbb{R}$ be as before in Definition \ref{plays-def}. The Lemma
  follows from the following calculation:
  \begin{eqnarray*}
    {\rm Col}(\aleph_1,\lambda)/\mbb{P} & \cong & (\mbb{R}\times{\rm
      Col}(\aleph_1,\lambda))/\mbb{P}\;\mbox{ (by Lemma
      \ref{coll-iso})}\\
    & \cong & (\mbb{R}/\mbb{P})\times{\rm Col}(\aleph_1,\lambda)
  \end{eqnarray*}
  and this last product is $\sigma$-closed by Lemma \ref{finsegsigma}.
\end{proof}

Lemma \ref{quot-sigma} points out the crucial difference between
strongly and weakly game-closed forcings: the quotient ${\rm
  Col}(\aleph_1,\lambda)/\mbb{P}$ will generally not be
$\sigma$-closed if $\mbb{P}$ is only weakly
$(\omega_1+1)$-game-closed.

Finally, we introduce the two remaining notions listed at the
beginning of the section. For an infinite cardinal $\kappa$, $\mbb{P}$
is called {\em $\kappa$-closed} if any $\mbb{P}$-descending chain of
length less than $\kappa$ has a lower bound in $\mbb{P}$. The poset
$\mbb{P}$ is called {\em $\kappa$-directed-closed} if it is closed
under directed subsets of size less than $\kappa$.
\cite{koenig04:_fragm_maxim} proves that {\rm PFA} is preserved by
$\omega_2$-closed forcings and \cite{larson00:_separ} that {\rm MM} is
preserved by $\omega_2$-directed-closed forcings.

\section{Indestructibility of set-theoretic
 axioms}\label{sec:indestruct}
\subsection{Forcing axioms}\label{subsec:forcingax}

This section should be seen as a continuation of work started in
\cite{koenig04:_fragm_maxim}. There it was shown that ${\rm PFA}$
implies failure of the approachability property at $\aleph_1$. While
we do not repeat the full proof here, it is interesting to mention
that ${\rm PFA}$ is applied to a poset $\mbb{Q}_0*\mbb{Q}_1*\mbb{Q}_2$
in this argument, where $\mbb{Q}_0$ adds a Cohen real, $\mbb{Q}_1$ is
a collapse with countable conditions, and $\mbb{Q}_2$ is specializing
a tree of size $\aleph_1$. The curious fact about the proof is that it
seems necessary for technical reasons to add the Cohen real right at
the start. In the following, we want to argue that the Cohen real is
absolutely necessary. We introduce the notion of a {\em
  $\Sigma$-collapsing} poset $\mbb{Q}$ which means that $\mbb{Q}$ can
be split up into two parts, where the first collapses everything in
sight without adding countable sequences, while the second is
arbitrary. Almost all known applications of ${\rm PFA}$ or ${\rm MM}$
are using $\Sigma$-collapsing posets. This section wants to point out
the few arguments where the forcing axiom for $\Sigma$-collapsing
posets is not enough even though the full forcing axiom suffices. In
other words, in the presented examples it is absolutely necessary to
add countable sequences before collapsing everything to $\aleph_1$.

\begin{definition}
  A poset $\mbb{Q}$ is called {\em $\Sigma$-collapsing} if it is of
  the form $\mbb{Q}=\mbb{Q}_0*\mbb{Q}_1$, where
  \begin{enumerate}
  \item $\mbb{Q}_0$ is $\aleph_1$-distributive and
  \item $\Vdash_{\mbb{Q}_0}|\mbb{Q}_1|\leq\aleph_1$.
  \end{enumerate}
  The class of $\Sigma$-collapsing posets that preserve stationary
  subsets of $\omega_1$ is denoted by $\GSC$.
\end{definition}

The point of this definition is that a $\Sigma$-collapsing poset will
typically collapse its own cardinality to $\aleph_1$ without adding
countable sequences. The final segment $\mbb{Q}_1$ is allowed to be
anything of size at most $\aleph_1$ though. Notice that, in the
context of forcing axioms, the class $\GSC$ in particular contains all
posets that are
\begin{itemize}
\item $\aleph_1$-distributive (take $\mbb{Q}_1$ to be trivial).
\item proper not adding reals (this implies
  $\aleph_1$-distributivity).
\item $\sigma\mbox{-closed}\,*\,{\rm ccc}$ (since we can assume the
  ${\rm ccc}$-poset to be of size at most $\aleph_1$
  \cite[p.62]{Kunen-intro}).
\end{itemize}
For example, ${\rm PFA}(\Gamma_\Sigma)$ implies the axiom ${\rm MRP}$
from \cite{moore05:_set} and ${\rm MM}(\Gamma_\Sigma)$ implies ${\rm
  SRP}$ (see e.g. \cite{MR1713438}). As mentioned in the introduction,
many classical applications of ${\rm PFA}$ are actually consequences
of the forcing axiom for posets of the form
$\sigma\mbox{-closed}\,*\,{\rm ccc}$ and all these are also included
in ${\rm PFA}(\Gamma_\Sigma)$.  Before we prove the main theorem of
this section, we reproduce the following Lemma. The proof is actually
a nice exercise but can also be looked up in \cite{kunen78:_satur}.

\begin{lemma}\label{kunen-lem}
  Let $\lambda$ be regular uncountable. Assume that $\mbb{R}$ is
  $\lambda$-closed, $A$ a relation on $\lambda$ and $\varphi$ a
  $\Sigma^1_1$-sentence. If there is a condition $r\in\mbb{R}$ such
  that $r\Vdash_\mbb{R}"(\lambda,A)\models\varphi"$ then
  $(\lambda,A)\models\varphi$. \hfill\qed
\end{lemma}

Let us prove the main preservation result. We show that strongly
$(\omega_1+1)$-game-closed forcings preserve the fragment of ${\rm
  MM}$ that contains all the $\Sigma$-collapsing posets. This will
later be used for interesting new independence results.

\begin{thm}\label{pres-semiproper}
  Assume that $V\models {\rm MM}(\GSC)$ and $\mbb{P}$ is strongly
  $(\omega_1+1)$-game-closed. Then $V^\mbb{P}\models{\rm MM}(\GSC)$.
\end{thm}
\begin{proof}
  Assume that
  \begin{equation}
    \label{eq:pres-semiproper0}
    \Vdash_\mbb{P}"\mbb{Q}\mbox{ is
      $\Sigma$-collapsing and preserves stationary subsets of $\omega_1$,"}
  \end{equation}
  where $\mbb{Q}_0*\mbb{Q}_1$ witnesses that $\mbb{Q}$ is
  $\Sigma$-collapsing and $\mbb{Q}$ is a $\mbb{P}$-name for a partial
  ordering. Let $\dot{\tau}_\xi\;(\xi<\omega_1)$ be a sequence of
  $\mbb{P}$-names for dense subsets of $\mbb{Q}$ and define the dense
  subsets of $\mbb{P}*\mbb{Q}$:
  \begin{equation}
    \label{eq:pres-semiproper1}
    D_\xi=\{(p,q):p\Vdash_\mbb{P}q\in\dot{\tau}_\xi\}.
  \end{equation}
  Now remember from the definitions that
  \begin{equation}
    \label{eq:pres-semiproper2}
    \Vdash_{\mbb{P}*\mbb{Q}_0}|\mbb{Q}_1|\leq\aleph_1
  \end{equation}
  and recall the poset $\mbb{R}$ from Definition \ref{plays-def} which
  is induced by $\mbb{P}$. It was shown in Lemma \ref{finsegsigma}
  that $\mbb{R}/\mbb{P}$ is $\sigma$-closed.
  \begin{claim}\label{stat-pres-claim}
    $\mbb{Q}_1$, as a forcing notion in $V^{\mbb{P}*\mbb{Q}_0}$,
    preserves all stationary subsets of $\omega_1$ in $V^\mbb{P}$.
  \end{claim}
  \begin{proof}
    This simply follows from the fact that
    $\mbb{Q}=\mbb{Q}_0*\mbb{Q}_1$ preserves all stationary subsets of
    $\omega_1$ in $V^\mbb{P}$ by (\ref{eq:pres-semiproper0}).
  \end{proof}
  \begin{claim}\label{four-step}
    The iteration $\mbb{P}*\mbb{Q}_0*\mbb{R}/\mbb{P}*\mbb{Q}_1$
    preserves stationary subsets of $\omega_1$.
  \end{claim}
  \begin{proof}
    Let $E\sub\omega_1$ be stationary. Clearly, the stationarity of
    $E$ is preserved in the three-step iteration
    $\mbb{P}*\mbb{Q}_0*\mbb{P}/\mbb{R}$. Now assume towards a
    contradiction that $\mbb{Q}_1$, as a forcing notion in
    $V^{\mbb{P}*\mbb{Q}_0*\mbb{R}/\mbb{P}}$, would destroy the
    stationarity of $E$. Then in $V^{\mbb{P}*\mbb{Q}_0}$
    \begin{equation}
      \label{eq:pres-semiproper3}
      \Vdash_{\mbb{R}/\mbb{P}}\mbox{ "there is a $\mbb{Q}_1$-name
        $\dot{C}$ for a club disjoint from $E$."}
    \end{equation}
    Using density arguments, it is straightforward to check that the
    quoted statement in (\ref{eq:pres-semiproper3}) is $\Sigma^1_1$
    over the structure
    $(\omega_1\cup\mbb{Q}_1,\in,E,\leq_{\mbb{Q}_1})$. So we can apply
    Lemma \ref{kunen-lem} and conclude that in $V^{\mbb{P}*\mbb{Q}_0}$
    \begin{equation}
      \label{eq:pres-semiproper4}
      \mbox{there is a $\mbb{Q}_1$-name $\dot{C}$ for a club
        disjoint from $E$.}
    \end{equation}
    But (\ref{eq:pres-semiproper4}) says that $\mbb{Q}_1$, as a
    forcing notion in $V^{\mbb{P}*\mbb{Q}_0}$, destroys the
    stationarity of $E$. This contradicts Claim \ref{stat-pres-claim}
    so we finished the proof of Claim \ref{four-step}.
  \end{proof}
  Now use the forcing axiom in the ground model to get a filter
  \begin{equation}
    \label{eq:pres-semiproper5}
    G*H_0*K*H_1\sub\mbb{P}*\mbb{Q}_0*\mbb{R}/\mbb{P}*\mbb{Q}_1
  \end{equation}
  that is sufficiently generic, in particular for all dense sets
  $D_\xi\;(\xi<\omega_1)$.
  \begin{claim}\label{G-ext-claim}
    The filter $G$ extends to a condition $q$ in $\mbb{P}$.
  \end{claim}
  \begin{proof}
    This is using the fact that $K$ yields a play of the Banach-Mazur
    game of length $\omega_1+1$ that exhausts $G$ in the sense that
    all elements of $G$ are refined during that play. But this play,
    given by $K$, follows Nonempty's winning strategy so there is a
    condition $q$ stronger than every condition in $G$.
  \end{proof}
  Claim \ref{G-ext-claim} finishes the proof since
  \begin{equation}
    \label{eq:pres-semiproper6}
    q\Vdash_\mbb{P}"H_0*H_1\sub\mbb{Q}_0*\mbb{Q}_1
    \mbox{ is generic for }\dot{\tau}_\xi\;(\xi<\omega_1)."
  \end{equation}
\end{proof}

To illustrate the significance of Theorem \ref{pres-semiproper}, we
turn to the following theorem of Magidor (see
\cite{cummings01:_squar}).

\begin{thm}\label{Magidor}
  ${\rm MM}$ implies the failure of ${\rm AP}_{\aleph_\omega}$.
  \hfill\qed
\end{thm}

In his proof, Magidor applies a forcing that shoots a new
$\omega$-sequence through $\aleph_\omega$. A Corollary of Theorem
\ref{pres-semiproper} explains why it is necessary in his argument to
add a new countable sequence high up in the set-theoretic hierarchy:

\begin{corollary}
  The following theory is consistent: $${\rm MM}(\aleph_\omega)+{\rm
    MM}(\GSC)+{\rm AP}_{\aleph_\omega}.$$
  \hfill\qed
\end{corollary}
\begin{proof}
  We add an ${\rm AP}_{\aleph_\omega}$-sequence to a model of ${\rm
    MM}$. The corollary now follows from Lemmas
  \ref{pres-bounded},\ref{force-AP} and Theorem \ref{pres-semiproper}.
\end{proof}

It has already been proved that ${\rm PFA} + {\rm AP}_{\aleph_\omega}$
is consistent. See \cite{cummings01:_squar} for more history on that.
Next we investigate the $\Sigma$-collapsing fragment of ${\rm PFA}$.

\begin{corollary}
  ${\rm Con}({\rm BPFA}+{\rm PFA}(\GSC)+{\rm AP}_{\aleph_1})$.
\end{corollary}
\begin{proof}
  By Lemmas \ref{pres-bounded},\ref{force-AP} and Theorem
  \ref{pres-semiproper} since $\mbb{A}_{\aleph_1}$ is strategically
  $\sigma$-closed and therefore proper.
\end{proof}

But PFA implies the failure of ${\rm AP}_{\aleph_1}$
\cite{koenig04:_fragm_maxim}, so we get:

\begin{corollary}\label{coro-reals}
  ${\rm BPFA}$ and $ {\rm PFA}(\GSC)$ together do not imply full ${\rm
    PFA}$. \hfill\qed
\end{corollary}

For transparency, we include a small chart that sums up the results in
this area. The class $\Gamma_{\rm cov}$ is the collection of all
posets that preserve stationary subsets of $\omega_1$ and have the
{\em covering property}, i.e. every countable set of ordinals in the
extension can be covered by a countable set in the ground model. Note
that $\GSC$ is a proper subset of $\Gamma_{\rm cov}$.

\medskip
\begin{tabular}{l|l}
  If $V\models{\rm MM}$ and $\mbb{P}$ is ... & then ... is true in
  $V^\mbb{P}$.\\
  \hline\\
  (1) $\omega_2$-distributive & ${\rm BMM}$\\
  (2) weakly $(\omega_1+1)$-game-closed & ${\rm BMM}$ + saturation of
  ${\rm NS}_{\omega_1}$\\
  (3) strongly $(\omega_1+1)$-game-closed & ${\rm BMM}$ + ${\rm
    MM}(\GSC)$\\
  (4) $\omega_2$-closed & ${\rm BMM}$ + ${\rm MM}(\Gamma_{\rm cov})$\\
  (5) $\omega_2$-directed-closed & ${\rm MM}$\\
\end{tabular}

\medskip In the table above, (1) is Lemma \ref{pres-bounded}, (2) is
in \cite{velickovic92:_forcin}, (3) is Theorem \ref{pres-semiproper},
(4) is in \cite{koenig04:_fragm_maxim}, and (5) is folklore but can be
looked up in \cite{larson00:_separ}. The papers
\cite{koenig04:_fragm_maxim} and \cite{koenig05:_kurep_namba} give
counterexamples which show that the results in this chart are
basically optimal. For example, it is shown in
\cite{koenig04:_fragm_maxim} that adding an ${\rm
  AP}_{\aleph_1}$-sequence with a strongly $(\omega_1+1)$-game-closed
forcing makes ${\rm PFA}$ (and therefore ${\rm MM}(\Gamma_{\rm cov})$)
fail in the extension. In \cite{koenig05:_kurep_namba}, an
$\omega_2$-closed forcing is constructed which makes ${\rm MM}$ fail
in the extension. More details on this last fact can actually be found
in Section \ref{sec:Kurepa-tree}\ref{subsec:reg-Ktrees} of this paper.

\subsection{Levy-Collapse}

There are similar preservation results for statements true in the
Levy-Collapse of a large cardinal. In \cite{koenig04:_gener}, an
attempt was made to axiomatize the theory of the Levy-Collapse with
the help of a reflection principle that is in the style of the
well-investigated reflection principles for stationary sets and is
actually a strengthening of these. This axiom, the {\em Game
  Reflection Principle} or ${\rm GRP}$ for short, proves all the
typical statements known to hold in the Levy-Collapse and is
formulated in terms of games. We need the following notions from
\cite{koenig04:_gener}:

\begin{definition}\label{gameondef}
Let $\theta$ be a regular cardinal.\hfill
\begin{enumerate}
\item (The Games) If $\mcal{A}\sub{^{<\omega_1}\theta}$ then the game
  $\mathbb{G}(\mcal{A})$\label{mathbbG(mathcalA)} has length
  $\omega_1$ and is played as follows: \medskip

\begin{tabular}{c|cccccc}
  I & $\alpha_0$ & $\alpha_1$ & $\ldots$ & $\alpha_\xi$ &
  $\alpha_{\xi+1}$ &
  $\ldots$\\\hline
  II & $\qquad \beta_0$ & $\qquad \beta_1$ & $\qquad\ldots$ &
  $\qquad \beta_\xi$ & $\qquad \beta_{\xi+1}$ & $\qquad\ldots$
\end{tabular}
\medskip\newline both players I and II play ordinals below $\theta$
and $$\mbox{II wins iff }\langle
\alpha_\xi,\beta_\xi:\xi<\omega_1\rangle\in[\mcal{A}],$$
where
$[\mcal{A}]=\{f\in{^{\omega_1}\theta}:f\restr
\xi\in\mcal{A}\mbox{ for all }\xi<\omega_1\}$.
\item (Restricted Games) For $B\sub\Hk{\lambda}$, define the game
  $\mathbb{G}^B(\mcal{A})$ by letting the winning conditions be the
  same as in $\mathbb{G}(\mcal{A})$ but imposing the restriction on
  both players to play ordinals in $B\cap\theta$.
\item ($\epsilon$-Approachability) A substructure $M\prec\Hk{\lambda}$
  of size $\aleph_1$ is called {\em $\epsilon$-approach\-able} if it
  is the limit of an $\epsilon$-chain of countable elementary
  substructures, i.e. there is an $\epsilon$-chain $\langle
  M_\xi:\xi<\omega_1\rangle$ with $M=\bigcup_{\xi<\omega_1}M_\xi$. We
  denote the set of all $\epsilon$-approachable substructures of
  $\Hk{\lambda}$ of size $\aleph_1$ by ${\rm EA}_\lambda$ and we drop
  the subscript if it is clear from the context.
\end{enumerate}
\end{definition}

\begin{definition}
  The {\em Game Reflection Principle} or ${\rm GRP}$ is the following
  statement:
 
  \medskip\hspace{0.7cm}\parbox{10cm}{Let $\theta$ be regular,
    $\mcal{A}\sub{^{<\omega_1}\theta}$ and $\lambda$ much larger
    than $\theta$. If II has a winning strategy in the game
    $\mathbb{G}^M(\mcal{A})$ for every $M\in{\rm EA}_\lambda$, then
    II has a winning strategy in $\mathbb{G}(\mcal{A})$.}
\end{definition}

We have the following two theorems from \cite{koenig04:_gener}:

\begin{thm}\label{GRP-gsc}
  ${\rm GRP}$ is equivalent to saying that $\omega_2$ is generically
  supercompact by $\sigma$-closed forcing, i.e. for every regular
  $\lambda$ there is $\mathbb{P}\in\Gamma$ such that $V^{\mathbb{P}}$
  supports $j:V\arr M$ with ${\rm cp}(j)=\omega_2$,
  $j(\omega_2)>\lambda$, and $j''\lambda\in M$.\hfill\qed
\end{thm}

\begin{thm}\label{Coll-GRP}
  Assume that $\kappa$ is supercompact. Then $$V^{{\rm
      Coll}(\omega_1,<\kappa)}\models{\rm GRP}.$$
  \hfill\qed
\end{thm}

Theorem \ref{Coll-GRP} is really contained in the stronger Theorem
\ref{pres-GRP}, so we postpone the proof. The following is the
preservation argument:

\begin{thm}\label{pres-GRP}
  Assume that $\kappa$ is supercompact and $\dot{\mbb{P}}$ is a ${\rm
    Coll}(\omega_1,<\kappa)$-name for a strongly
  $(\omega_1+1)$-game-closed partial ordering. Then $$V^{{\rm
      Coll}(\omega_1,<\kappa)*\dot{\mbb{P}}}\models{\rm GRP}.$$
\end{thm}
\begin{proof}
  By Theorem \ref{GRP-gsc} it suffices to show that $\omega_2$ is
  generically supercompact by $\sigma$-closed forcing. For any regular
  $\lambda$ fix $j:V\arr M$ such that ${\rm cp}(j)=\kappa$,
  $j(\kappa)>\lambda$, and $j''\lambda\in M$. Without restriction,
  $\lambda>|\dot{\mbb{P}}|$.
  \begin{claim}\label{claim-sigma}
    $\Vdash_{{\rm Coll}(\omega_1,<\kappa)}{\rm
      Coll}(\omega_1,[\kappa,j\kappa))/\dot{\mbb{P}}$ is
    $\sigma$-closed.
  \end{claim}
  \begin{proof}
    Note first that this makes sense because
    $j\kappa>\lambda>|\dot{\mbb{P}}|$ and therefore it is forced by
    ${\rm Coll}(\omega_1,<\kappa)$ that $\dot{\mbb{P}}$ is a complete
    subalgebra of ${\rm Coll}(\omega_1,[\kappa,j\kappa))$. This last
    statement is basically Lemma \ref{str-embeds}. The claim now
    follows from Lemma \ref{quot-sigma}.
  \end{proof}
  Now set $\mcal{Q}={\rm Coll}(\omega_1,<\kappa)*\dot{\mbb{P}}$.
  Standard arguments yield
  \begin{eqnarray*}
    j\mcal{Q} & = & j({\rm Coll}(\omega_1,<\kappa)*\dot{\mbb{P}})\\
    & = & {\rm Coll}(\omega_1,<j\kappa)*j\dot{\mbb{P}}\\
    & = & {\rm Coll}(\omega_1,<\kappa)*{\rm
      Coll}(\omega_1,[\kappa,j\kappa))*j\dot{\mbb{P}}\\
    & = & {\rm Coll}(\omega_1,<\kappa)*\dot{\mbb{P}}*{\rm
      Coll}(\omega_1,[\kappa,j\kappa))/\dot{\mbb{P}}*j\dot{\mbb{P}}\\
    & = & \mcal{Q}*{\rm
      Coll}(\omega_1,[\kappa,j\kappa))/\dot{\mbb{P}}*j\dot{\mbb{P}}\\
    & = & \mcal{Q}*j\mcal{Q}/\mcal{Q}.
  \end{eqnarray*}
  By the Extension Lemma \ref{silverlem} we can, in $V^{j\mcal{Q}}$,
  extend $j$ to $$j^*:V^\mcal{Q}\arr M^{j\mcal{Q}}.$$
  Notice first
  that $j''\lambda\in M^{j\mcal{Q}}$. Finally, it holds in the model
  $V^\mcal{Q}$ that $\omega_2$ is generically $\lambda$-supercompact
  by the forcing $j\mcal{Q}/\mcal{Q}$. But remember that
  $$j\mcal{Q}/\mcal{Q}={\rm
    Coll}(\omega_1,[\kappa,j\kappa))/\dot{\mbb{P}}*j\dot{\mbb{P}}$$
  is
  $\sigma$-closed by Claim \ref{claim-sigma}.
\end{proof}

Note that ${\rm GRP}$ implies ${\rm CH}$ \cite{koenig04:_gener} and
therefore ${\rm AP}_{\aleph_1}$. So Theorem \ref{pres-GRP} tells us
nothing if we take $\dot{\mbb{P}}$ to be the usual forcing that adds
an ${\rm AP}_{\aleph_1}$-sequence. But we can get the following
interesting corollary:

\begin{corollary}\label{APomegaGRP}
  ${\rm Con}({\rm GRP} + {\rm AP}_{\aleph_\omega})$.
\end{corollary}
\begin{proof}
  By Lemma \ref{force-AP} and Theorem \ref{pres-GRP}.
\end{proof}

Corollary \ref{APomegaGRP} is interesting because Shelah
\cite{shelah79} has shown that if $\kappa$ is supercompact then ${\rm
  AP}_{\kappa^{+\omega}}$ fails. This contrasts the above result in
the sense that generic supercompactness of $\omega_2$ by
$\sigma$-closed forcing does not imply that ${\rm AP}_{\aleph_\omega}$
fails. To end the section, let us remark that weakly
$(\omega_1+1)$-game-closed forcings can introduce square-sequences and
this would violate ${\rm GRP}$ in a very strong fashion
\cite{koenig04:_gener}. So we cannot hope to have a preservation
theorem for weakly $(\omega_1+1)$-game-closed forcings.

\section{Kurepa-trees}\label{sec:Kurepa-tree}

\subsection{Jonsson cardinals}\label{subsec:Jonsson}

Let us go back to Theorem \ref{Magidor} for a while. Magidor's
argument that ${\rm MM}$ implies failure of the approachability
property at $\aleph_\omega$ has raised hopes that ${\rm MM}$ can
provide a good picture of the combinatorics of the cardinal
$\aleph_\omega$. Even earlier \cite{foreman88}, it was shown that
${\rm MM}$ implies ${\rm SCH}$ which puts severe restrictions on the
cardinal arithmetic at $\aleph_\omega$. In the light of all this, it
seemed possible that ${\rm MM}$ implies $\aleph_\omega$ is Jonsson
which would solve an old question. We give an argument to refute this
last implication, i.e. we show that ${\rm MM}$ is consistent with
${\aleph_\omega}$ not being Jonsson. This answers a question raised
during the Oberwolfach set theory meeting in 2005. The following
theorem is well-known, see for example \cite{higherinf}.

\begin{thm}\label{least-Jonsson}
  Assume that $\mu$ is the least Jonsson cardinal. Then $\mu$ is
  $\lambda$-Rowbottom for some $\lambda<\mu$.\hfill\qed
\end{thm}

Kurepa-trees are natural counterexamples to model-theoretic transfer
properties. The next lemma is probably standard, but we give the proof
for convenience.

\begin{lemma}\label{notcc-lem}
  Let $\kappa<\lambda<\mu$ be cardinals where $\lambda$ is regular and
  assume that there is a $\lambda$-Kurepa-tree $T$ with at least
  $\mu$-many cofinal branches. Then
  $(\mu,\lambda)\notcc(\lambda,\kappa)$.
\end{lemma}
\begin{proof}
  Let $\mcal{B}$ be a collection of $\mu$-many cofinal branches of
  $T$. The structure $(\mcal{B},T)$ is of type $(\mu,\lambda)$, so
  suppose towards a contradiction that there exists
  $$(\mcal{A},S)\prec(\mcal{B},T),$$
  where $\mcal{A}\sub\mcal{B}$ is
  of size $\lambda$ and $S\sub T$ is of size $\kappa$. Find
  $\delta<\lambda$ such that $S\sub T_{<\delta}$. Then by
  elementarity, every two branches in $\mcal{A}$ split within the
  structure $(\mcal{A},S)$ which implies that $T_\delta$ would have
  size at least $|\mcal{A}|=\lambda$. This contradicts the fact that
  levels of $T$ have size less than $\lambda$.
\end{proof}

Let us point out again (cf. Section
\ref{sec:indestruct}\ref{subsec:forcingax}) that ${\rm MM}$ is
preserved by $\omega_2$-directed-closed forcings. This is used
crucially in the proof of the next theorem.

\begin{thm}\label{MM-Jonsson}
  ${\rm MM}$ does not imply that $\aleph_\omega$ is Jonsson.
\end{thm}
\begin{proof}
  We construct a model of ${\rm MM}$ in which $\aleph_\omega$ is not
  Jonsson. First note that ${\rm MM}$ is consistent with
  $2^\lambda=\lambda^+$ for all $\lambda\geq\aleph_1$ since this
  instance of the ${\rm GCH}$ can be forced with an
  $\omega_2$-directed-closed forcing. So we start with a model
  $$V\models{\rm MM}+2^\lambda=\lambda^+\mbox{ for all
  }\lambda\geq\aleph_1$$
  and define a full support Easton product
  $$\mbb{K}_\omega=\prod_{1<n<\omega}\mbb{K}_n,$$
  where $\mbb{K}_n$ is
  the usual forcing to add an $\aleph_n$-Kurepa-tree with
  $\aleph_\omega$-many branches (see e.g. \cite{Kunen-intro}). Note
  that $\mbb{K}_n$ is $\aleph_n$-directed-closed and has the
  $\aleph_{n+1}$-chain condition as we assumed the arithmetic
  $2^{<\aleph_n}=\aleph_n$. The usual arguments for the Easton product
  yield that $\mbb{K}_\omega$ preserves all cardinals
  $\aleph_n\;(n<\omega)$ and therefore preserves $\aleph_\omega$. It
  is also easy to see that $\mbb{K}_\omega$ is an
  $\omega_2$-directed-closed forcing and hence preserves ${\rm MM}$.
  By Lemma \ref{notcc-lem} we have
  \begin{equation}
    \label{eq:MM-Jonsson1}
    (\aleph_\omega,\aleph_m)\notcc(\aleph_m,\aleph_n)\mbox{ for
      all }n<m<\omega
  \end{equation}
  in the generic extension $V^{\mbb{K}_\omega}$. Now assume that
  $\aleph_\omega$ is Jonsson in $V^{\mbb{K}_\omega}$ and use Theorem
  \ref{least-Jonsson}. We get that
  $(\aleph_\omega,\aleph_m)\cc(\aleph_\omega,\aleph_n)$ holds for some
  $n<m<\omega$. This contradicts (\ref{eq:MM-Jonsson1}).
\end{proof}

The proof of Theorem \ref{MM-Jonsson} actually shows that ${\rm
  MM}^{++}$ does not imply that $\aleph_\omega$ is Jonsson, where
${\rm MM}^{++}$ means that $\omega_1$-many names for stationary
subsets of $\omega_1$ can be pulled back into the ground model. We
generally tried to avoid these 'plus-versions' of forcing axioms, the
interested reader is referred to
\cite{baumgartner84:_applic_proper_forcin_axiom} or \cite{MR1713438}.

\subsection{Regressive Kurepa-trees}
\label{subsec:reg-Ktrees}

The notion of a regressive Kurepa-tree was introduced in
\cite{koenig05:_kurep_namba} in order to answer the question if ${\rm
  MM}$ is sensitive to $\omega_2$-closed forcings. Surprisingly, ${\rm
  MM}$ turned out to be sensitive to $\omega_2$-closed forcings but
only the Namba-fragment of ${\rm MM}$ can be violated. The key notion
was that of an $\omega_1$-regressive $\omega_2$-Kurepa-tree which can
be added by an $\omega_2$-closed forcing and it was shown that ${\rm
  MM}$ is false once such a tree is added. We want to point out in
this section that regressive Kurepa-trees have strong impact on higher
versions of Chang's Conjecture even though they can be added with a
very mild forcing.

\begin{definition}
  For any tree $T$ say that the level $T_\alpha$ is {\em
    non-stationary} if there is a function $f_\alpha:T_\alpha\arr
  T_{<\alpha}$ which is {\em regressive} in the sense that
  $f_\alpha(x)<_Tx$ for all $x\in T_\alpha$ and if $x,y\in T_\alpha$
  are distinct then $f_\alpha(x)$ or $f_\alpha(y)$ is strictly above
  the meet of $x$ and $y$.
\end{definition}

\begin{definition}
  Let $X$ be a set of ordinals. A $\lambda$-Kurepa-tree $T$ will be
  called {\em $X$-regressive} if for all limit ordinals
  $\alpha<\lambda$ with ${\rm cf}(\alpha)\in X$ the level $T_\alpha$
  is non-stationary.
\end{definition}

The following two theorems appear in \cite{koenig05:_kurep_namba}.

\begin{thm}\label{add-regKtree1}
  For all uncountable regular $\lambda$ there is a $\lambda$-closed
  forcing that adds a $\lambda$-regressive $\lambda$-Kurepa-tree.
  \hfill\qed
\end{thm}

\begin{thm}\label{MM-no-reg-Ktree}
  Under {\rm MM}, there are no $\omega_1$-regressive
  $\lambda$-Kurepa-trees for any uncountable regular $\lambda$.
  \hfill\qed
\end{thm}

A close examination of the proof of Theorem \ref{add-regKtree1}
actually gives:

\begin{thm}\label{addregKtree2}
  Assume $2^{<\lambda}=\lambda$ and $\kappa<\lambda<\mu$, where
  $\kappa$ and $\lambda$ are regular. There is a
  $\kappa^+$-directed-closed, $\lambda$-closed,
  $\lambda^+-\mathit{cc}$ forcing that adds a
  $[\kappa,\lambda)$-regressive $\lambda$-Kurepa-tree with at least
  $\mu$-many branches. \hfill\qed
\end{thm}

Regressive Kurepa-trees are even stronger counterexamples to
model-theoretic transfer properties than the regular Kurepa-trees
considered in Section \ref{sec:Kurepa-tree}\ref{subsec:Jonsson}. This
can be seen from the following Lemma.

\begin{lemma}\label{regKtree-cc}
  Let $\kappa<\lambda$ be regular. Assume there is a
  $\{\kappa\}$-regressive $\lambda$-Kurepa-tree $T$ with at least
  $\mu$-many branches and suppose $\kappa^{<\kappa}=\kappa$. Then
  $(\mu,\lambda)\notcc(\kappa^+,\kappa)$.
\end{lemma}
\begin{proof}
  Let $\mcal{B}$ be the set of cofinal branches of $T$ and consider
  the structure $(\mcal{B},T)$ which is of type $(\mu,\lambda)$. Now
  assume towards a contradiction that
  $(\mu,\lambda)\cc(\kappa^+,\kappa)$ would hold, so we find a
  substructure $$(\mcal{A},S)\prec(\mcal{B},T),$$
  where $\mcal{A}$ has
  size $\kappa^+$ and $S$ has size $\kappa$. Define $\delta=\sup({\rm
    ht}"S)$, we have two cases:

  {\bf Case 1:} if ${\rm cf}(\delta)=\kappa$ then $T_\delta$ is a
  non-stationary level of the tree $T$. A straightforward argument
  using the fact that there is a regressive 1-1 function defined on
  $T_\delta$ shows that $\mcal{A}$ has size at most $|S|=\kappa$. This
  is a contradiction.

  {\bf Case 2:} if ${\rm cf}(\delta)<\kappa$ then $S$ has a cofinal
  subtree $S_0$ of height $\kappa_0<\kappa$. Since $|S|=\kappa$, the
  number of branches through $S_0$ can not be larger than
  $\kappa^{\kappa_0}\leq\kappa^{<\kappa}=\kappa$. Again,
  contradiction.
\end{proof}

\begin{corollary}\label{MM-cc}
  MM is consistent with
  $$(\aleph_{m+1},\aleph_m)\notcc(\aleph_{n+1},\aleph_n)$$
  for all $1<n<m$ simultaneously.
\end{corollary}
\begin{proof}
  Using an Easton product similar to the proof of Theorem
  \ref{MM-Jonsson}: we start with a model of $"{\rm
    MM}+2^\lambda=\lambda^+$ for all $\lambda\geq\aleph_1$". Then for
  all $m>2$ we add an $[\aleph_2,\aleph_m)$-regressive
  $\aleph_m$-Kurepa-tree with $\aleph_{m+1}$-many branches. This
  product is $\omega_2$-directed-closed by Theorem \ref{addregKtree2}.
  Notice that in the extension we have $2^{<\aleph_n}=\aleph_n$ for
  all $1<n<\omega$. The statement of the corollary then follows from
  Lemma \ref{regKtree-cc}.
\end{proof}

\bibliography{locco}\bibliographystyle{plain}

\end{document}